\documentclass[10pt]{amsart}
\usepackage{amsmath,amsfonts,amssymb}
\usepackage[english,francais]{babel}
\newtheorem{theorem}{Theorem}

\newtheorem{e-proposition}[theorem]{Proposition}
\newtheorem{corollary}[theorem]{Corollary}
\newtheorem{e-definition}[theorem]{Definition\rm}
\newtheorem{remark}{\it Remark\/}

\def\weak{weak-* }
\def\teta{\vartheta}
\def\fhi{\varphi}
\def\epsi{\varepsilon}
\def\ti{-\allowhyphens}
\def\lra{\longrightarrow}
\def\hb{{\rm H}_{\rm b}}
\def\one{\mathbf{1\kern-1.6mm 1}}
\def\No{N\raise4pt\hbox{\tiny o}\kern+.2em}
\def\no{n\raise4pt\hbox{\tiny o}\kern+.2em}

\def\beq{\begin{equation}}
\def\eeq{\end{equation}}
\def\deq{\stackrel{\mbox{\tiny def}}{=}}
\begin{document}
\selectlanguage{english}
\title[co-amenability for groups and von Neumann algebras]{On co-amenability for groups\\ and von Neumann algebras}
\author{Nicolas Monod \& Sorin Popa}
\address{N.M.: University of Chicago\\ Chicago IL 60637 USA}
\email{monod@math.uchicago.edu}
\address{S.P.: UCLA\\ Los Angeles CA 90095 USA}
\email{popa@math.ucla.edu}
\begin{abstract}
We first show that co-amenability does not pass to subgroups, answering a question asked by Eymard in~1972. We then address co-amenability for von Neumann algebras, describing notably how it relates to the former.

\vskip2mm\noindent\selectlanguage{french}
\scshape R\'esum\'e.\upshape~
Nous d\'emontrons tout d'abord que la comoyennabilit\'e ne passe pas aux sous-groupes, r\'esolvant ainsi une question pos\'ee par P.~Eymard en~1972. Nous \'etudions ensuite la comoyennabilit\'e pour les alg\`ebres de von Neumann; en particulier, nous clarifions le lien avec le cas des groupes.
\selectlanguage{english}
\end{abstract}
\maketitle

\selectlanguage{english}
\noindent
\textbf{\Large Co-amenable subgroups}

A subgroup $H$ of a group $G$ is called \emph{co-amenable} in $G$ if it has the following relative fixed point property:

\slshape Every continuous affine $G$\ti action on a convex compact subset of a locally convex space with an $H$\ti fixed point has a $G$\ti fixed point.\upshape

\smallskip

This is equivalent to the existence of a $G$\ti invariant mean on the space $\ell^\infty(G/H)$ or to the weak containment of the trivial representation in $\ell^2(G/H)$; see~\cite{Eymard72} (compare also Proposition~\ref{prop_coh} below). Alternative terminology is \emph{amenable pair} or \emph{co-F\o lner}, and more generally Greenleaf~\cite{Greenleaf69} introduced the notion of an \emph{amenable action} (conflicting with more recent terminology) further generalized by Zimmer~\cite{Zimmer78} to amenable pairs of actions.

\smallskip

In the particular case of a normal subgroup $H\lhd G$, co-amenability is equivalent to amenability of the quotient $G/H$.

\medskip

One checks that for a triple
$$K<H<G,$$
co-amenability of $K$ in $H$ and of $H$ in $G$ implies that $K$ is co-amenable in $G$. In~1972, Eymard proposed the following problem~\cite[\No 3\S8]{Eymard72}:

\smallskip

\slshape If we assume $K$ co-amenable in $G$, is it also co-amenable in $H$?\upshape

\smallskip

(It is plain that $H$ is co-amenable in $G$.)

\medskip

If $K$ is normal in $G$, then the answer is positive, since it amounts to the fact that amenability passes to subgroups. We give an elementary family of counter-examples for the general case. Interestingly, we can even get some normality:

\begin{theorem}
\label{thm_groups}\ 

\noindent
(i)~There are triples of groups $K \lhd H \lhd G$ with $K$ co-amenable in $G$ but not in $H$.

\smallskip
\noindent
(ii)~For \textbf{any} group $Q$ there is a triple $K\lhd H_0 < G$ with $K$ co-amenable in $G$ and $H_0/K\cong Q$.

Moreover, all groups in~(i) can be assumed ICC (infinite conjugacy classes). In~(ii), the group $G$ is finitely generated/presented, ICC, torsion-free, etc. as soon as $Q$ has the corresponding property.
\end{theorem}

\begin{remark}
\label{rem_fg}%
We shall see that one can also find counter-examples $K<H_0<G$ with all three groups finitely generated.
\end{remark}

\begin{remark}We mention that in a recent independent work V.~Pestov also solved Eymard's problem (see V.~Pestov's note in this same issue).
\end{remark}

We start with the following general setting: Let $K$ be any group and $\teta:K\to K$ a monomorphism. Denote by $G=K*_\teta$ the corresponding HNN-extension, \emph{i.e.}
$$G\deq \big\langle K,t\ |\ t k t^{-1}=\teta(k)\ \forall\,k\in K\big\rangle.$$
Write $H\deq \bigcup_{k\in {\bf Z}} t^{-k}Kt^k$ so that $G=H\rtimes <t>$.

\begin{e-proposition}
\label{prop_HNN}%
The group $K$ is co-amenable in $G$.
\end{e-proposition}

\begin{proof}
Note that $H$ is co-amenable in $G$, so that by applying the relative fixed point property to the space of means on $\ell^\infty(G/K)$ we deduce that it is enough to find an $H$\ti invariant mean
$$\mu:\ \ell^\infty(G/K)\lra \bf{R}.$$
For every $k\in\bf{N}$ we define a mean $\mu_k$ by $\mu_k(f)=f(t^{-k}K)$ for $f$ in $\ell^\infty(G/K)$. This mean is left invariant by the group $t^{-k}Kt^k$. Since $H$ is the increasing union of $t^{-k}Kt^k$ for $k\to \infty$, any \weak limit of the sequence $\{\mu_k\}$ is an $H$\ti invariant mean.
\end{proof}

We set $H_0=t^{-1}Kt$ and consider $K<H_0<H\lhd G$. Observe that $K$ is not co-amenable in $H_0$ unless $\teta(K)$ is co-amenable in $K$, and hence also not in $H$ when $K\lhd H$. Thus we will construct our counter-examples by finding appropriate pairs $(K,\teta)$ with $\teta(K)$ not co-amenable in $K$:

\emph{For Theorem~\ref{thm_groups}}, let $Q$ be any group, let $K=\bigoplus_{n\geq 0} Q$ be the direct limit of finite Cartesian powers and let $\teta$ be the shift; then $\teta(K)$ is co-amenable in $K$ if and only if $Q$ is amenable. The claims follow; observe that $G$ is the restricted wreath product of $Q$ by $\bf{Z}$. \emph{For Remark~\ref{rem_fg}}, we may for instance take $K$ to be the free group on two generators and define $\teta$ by choosing two independent elements in the kernel of a surjection of $K$ to a non-amenable group.

\bigskip

Amenability of groups can be characterized in terms of \emph{bounded cohomology}; see~\cite{Johnson}. 
Here is a generalization, which shows that actually the transitivity and (lack of) hereditary properties of co-amenability were to be expected.

\begin{e-proposition}
\label{prop_coh}%
The following are equivalent:

\begin{itemize}
\item[(i)] $H$ is co-amenable in $G$.

\item[(ii)] For any dual Banach $G$\ti module $V$ and any $n\geq 0$, the restriction $\hb^n(G,V)\to \hb^n(H,V)$ is injective.
\end{itemize}
\end{e-proposition}

\begin{proof}
``$(i)\Rightarrow(ii)$''~is obtained by constructing a \emph{transfer} map (at the cochain level) using an invariant mean on $G/H$, along the lines of~\cite[\S8.6]{Monod}.

``$(ii)\Rightarrow(i)$''~Let $V$ be the $G$\ti module $\big(\ell^\infty(G/H)/{\bf C}\big)^*$, where ${\bf C}$ is embedded as constants. The extension $0\to V \to \ell^\infty(G/H)^* \to {\bf C} \to 0$ splits over $H$ (by evaluation on the trivial coset) so that in the associated long exact sequence  the transgression map ${\bf C}\to\hb^1(H,V)$ vanishes. By~(ii) and naturality, the transgression ${\bf C}\to\hb^1(G,V)$ also vanishes, and so the long exact sequence for $G$ yields a $G$\ti invariant element $\mu$ in $\ell^\infty(G/H)^*$ which does not vanish on the constants. Taking the absolute value of $\mu$ (for the canonical Riesz space structure) and normalizing it gives an invariant mean.
\end{proof}

\noindent
\textbf{\Large Co-amenable von Neumann subalgebras}

\nobreak
Let $N$ be a finite von Neumann algebra and $B\subset N$ a von Neumann subalgebra. For simplicity, we shall assume $N$ has countably decomposable centre. Equivalently, we assume $N$ has a normal faithful tracial state $\tau$. Let $B\subset N \stackrel{e_B}{\subset} \langle N,B\rangle$ be the basic construction for $B\subset N, \tau$, \emph{i.e.} $\langle N,B\rangle = JBJ' \cap \mathcal{B}\big(L^2(N,\tau)\big)$ with $J$ the canonical conjugation on $L^2(N,\tau)$ and $e_B$ the canonical projection of $L^2(N,\tau)$ onto $L^2(B, \tau)$. Note that, up to isomorphism, the inclusion $N \subset \langle N,B\rangle$ does not depend on $\tau$. We also denote by $\mathrm{Tr}$ the unique 
normal semifinite faithful trace on $\langle N,B\rangle$ given by 
$\mathrm{Tr}(xe_By)=\tau(xy), x, y \in N$. 

\begin{e-definition}[\cite{Popa86}, see also~\cite{Popa99}]
The subalgebra $B$ is \emph{co-amenable} in $N$ if there exists a norm one projection $\Psi$ of $\langle N,B\rangle$ onto $N$. One also says that $N$ is \emph{amenable relative to $B$}.
\end{e-definition}

\begin{remark}
\label{rem_alg}%
If $B_0\subset B \subset N$ and $N$ is amenable relative to $B_0$ then $N$ is amenable relative to $B$. Indeed, since $\langle N,B\rangle$ is included in $\langle N, B_0\rangle$, we may restrict $\Psi:\langle N,B_0\rangle\to N$ to $\langle N,B\rangle$.
\end{remark}

The next result from~\cite{Popa86,Popa99} provides several alternative characterizations of 
relative amenability. The proof, which we have included for 
the reader's convenience, is an adaptation to the case of inclusions of algebras  
of Connes' well known results for single von Neumann algebras~\cite{Connes76b}.

\begin{e-proposition}
\label{prop_conditions}%
With $B\subset N,\tau$ as before, the following conditions are equivalent:

\begin{itemize}
\item[(i)] $B$ is co-amenable in $N$.

\item[(ii)] There exists a state $\psi$ on $\langle N,B\rangle$ with $\psi(uXu^*) = \psi(X)$ for all $X\in \langle N,B\rangle$, $u\in\mathcal{U}(N)$.

\item[(iii)] For all $\epsi >0$ and all finite $F\subset \mathcal{U}(N)$ there is a projection $f$ in $\langle N,B\rangle$ with $\mathrm{Tr} f<\infty$ such that 
$\|ufu^* - f\|_{2,\mathrm{Tr}} < \epsi \|f\|_{2,\mathrm{Tr}}$ for all $u\in F$.
\end{itemize}
\end{e-proposition}

\noindent
Here~(ii) parallels the invariant mean criterion while~(iii) is a F\o lner type condition.

\begin{proof} ``$ (i) \Leftrightarrow (ii)$''~ 
If $\Psi$ is a norm one projection onto 
$N$, then by Tomiyama's Theorem it is a 
conditional expectation onto $N$; thus $\psi = \tau \circ \Psi$ 
is $\mathrm{Ad}(u)$-invariant $\forall u\in \mathcal{U}(N)$. Conversely, 
if $\psi$ satisfies $(ii)$ above and $X\in \langle N,B \rangle$ 
then let $\Psi(X)$ denote the unique element in $N$ satisfying 
$\tau(\Psi(X)x), \forall x\in N$. It is immediate to see 
that $X \mapsto \Psi(X)$ is a conditional expectation 
(thus a norm one projection) onto $N$.  

``$(ii) \Leftrightarrow (iii)$''~ Assuming $(iii)$, denote by $I$ the 
set of finite subsets of $\mathcal{U}(N)$ ordered by 
inclusion and for each $i\in I$ let $f_i \in \mathcal{P}(\langle N, B \rangle)$ 
be a non-zero projection with $\mathrm{Tr}(f_i) < \infty$ such that 
$\|uf_iu^* - f\|_{2,\mathrm{Tr}} < |i|^{-1} \|f\|_{2,\mathrm{Tr}}, \forall u\in i$. 
For 
each $X\in \langle N, B \rangle$ let 
$\psi(X) \deq \lim_i \mathrm{Tr}(Xf_i)/\mathrm{Tr}(f_i)$ for some Banach limit over $i$. Then 
$\psi$ is easily seen to be an $\mathrm{Ad}(u)$-invariant state on $\langle N,B \rangle$, 
$\forall u\in \mathcal{U}(N)$. 

Conversely, if $\psi$ satisfies $(ii)$ then by Day's convexity trick 
it follows that for any finite $F \subset \mathcal{U}(N)$ and any 
$\epsi > 0$ there exists $\eta \in L^1(\langle N, B\rangle, \mathrm{Tr})_+$ 
such that $\mathrm{Tr}(\eta)=1$ and $\|u\eta u^* - \eta\|_{1,\mathrm{Tr}} \leq \epsi$. 
By the Powers-St\o rmer inequality, this implies 
$\|u\eta^{1/2} u^* - \eta^{1/2}\|_{2,\mathrm{Tr}} \leq 
\epsi^{1/2} \|\eta^{1/2}\|_{2,\mathrm{Tr}}$. Thus, 
Connes' ``joint distribution'' trick (see~\cite{Connes76}) applies to get a spectral projection 
$f$ of $\eta^{1/2}$ such that $\|ufu^*-f\|_{2,\mathrm{Tr}} < \epsi^{1/4}\|f\|_{2,\mathrm{Tr}}$, 
$\forall u\in F$.  
\end{proof} 

The next result relates the co-amenability of subalgebras 
with the co-amenability of subgroups, generalizing 
a similar argument for single algebras/groups from~\cite{Connes76b}. 

\begin{e-proposition}
\label{prop_equi}%
Let $G$ be a (discrete) group, $(B_0, \tau_0)$ a finite 
von Neumann algebra with a normal faithful tracial state and $\sigma: G\to \mathrm{Aut}(B_0,\tau_0)$ a trace-preserving cocycle action of 
$G$ on $(B_0, \tau_0)$. Let $N=B_0\rtimes_{\sigma} G$ be the corresponding crossed product von Neumann algebra with its normal faithful tracial state given by $\tau(\sum_{g\in G}b_g u_g) = \tau_0(b_e)$. Let also 
$H< G$ be a subgroup and 
$B=B_0\rtimes_{\sigma} H$. 

\nobreak
Then $B$ is co-amenable in $N$ if and only if the group 
$H$ is co-amenable in $G$
\end{e-proposition}

In particular, taking $B_0={\bf C}$, we deduce:

\begin{corollary}
\label{cor_equi}%
Let $H<G$ be (discrete) groups. Then the inclusion of group von Neumann algebras $L(H)\subset L(G)$ 
is co-amenable if and only if $H$ is co-amenable in $G$.
\end{corollary}

\begin{proof}[Proof of Proposition~\ref{prop_equi}]
``$\Leftarrow$''~Let $F_i\nearrow G/H$ be a net of 
finite F\o lner sets, which we identify with 
some sets of representatives $F_i\subset G$. 
For each $i$ and $X\in \langle N,B\rangle$ set
$$\Psi_i(X) = |F_i|^{-1}\Phi\Big(\big(\sum_{g\in F_i}u_g e_B u_g^*)X(\sum_{h\in F_i}u_h e_B u_h^*\big)\Big),$$
where $\Phi: {\text{\rm sp}}Ne_BN \to N$ is the canonical operator valued 
weight given by $\Phi(xe_By) = xy$ for $x,y\in N$. Then clearly 
$\Psi_i$ are completely positive, normal, unital, $B_0$\ti bimodular 
and satisfy 
$\Psi_i(\sum_{g\in G} b_gu_g) = \sum_{g\in G}\alpha_g^i b_g u_g$, 
where $\alpha_g^i = |gF_i\cap F_i|/|F_i|$.   
Thus, if we put 
$\Psi(X)\deq \lim_i \Psi_i(X)$ for some Banach limit over $i$, 
then $\Psi(X) = X$ for all $X\in N$, $\Psi(\langle N,B\rangle)\subset N$ 
and $\Psi$ is completely positive and unital. Thus $\Psi$ is a 
norm one projection.

``$\Rightarrow$''~Let $\{g_i\}_i \subset G$ be a set of representatives 
for the classes $G/H$ and denote by 
$\mathcal{A}\subset \langle N,B\rangle$ the set of 
operators of the form $\sum_{i} \alpha_i u_{g_i} e_B u_{g_i}^*$, 
where $\alpha_i\in{\bf C}$, $\sup_{i}|\alpha_i|<\infty$. Clearly, $\mathcal{A}\cong\ell^\infty(G/H)$. Note that ${\rm Ad}(u_h)$ implements the left translation by $h\in G$ on $\ell^\infty(G/H)$ 
via this isomorphism. Put $\fhi \deq \tau\circ \Psi|_\mathcal{A}$; 
then $\fhi$ is clearly a left invariant mean on $\ell^\infty(G/H)$.
\end{proof}

\begin{corollary}
There are inclusions of type $\mathrm{II}_1$ factors $N_0\subset N\subset M$ such that $N_0$ is co-amenable in $M$ but not in $N$. In fact, there are such examples with $N$ having property~(T) relative to $N_0$ in the sense of~\cite{Anantharaman,Popa86,Popa99}.
\end{corollary}

\begin{proof}
Let $Q$ be an infinite group with property~(T) and let $K\lhd H_0 < G$ be as in point~(ii) of Theorem~\ref{thm_groups}. Set $N_0=L(K)$, $N=L(H_0)$, $M=L(G)$. Choosing $Q$ to be ICC ensures that these three algebras are factors of type $\mathrm{II}_1$. The claim follows now from Corollary~\ref{cor_equi}.
\end{proof}

We proceed now to present an alternative example of interesting 
co-amenable inclusions of type $\mathrm{II}_1$ factors.

Let $\omega\in\beta{\bf N}\setminus {\bf N}$ be a free ultrafilter on ${\bf N}$. Denote as usual by $R$ the hyperfinite $\mathrm{II}_1$ factor and let $R^\omega\deq \ell^\infty(N,R)/I_\omega$, where $I_\omega\deq \big\{(x_n)_{n\in{\bf N}} | \lim_{n\to\omega}\tau(x_n^*x_n)=0\big\}$. Then $R^\omega$ is well known to be a type $\mathrm{II}_1$ factor with its unique trace given by $\tau_\omega((x_n)_{n\in {\bf N}}) = \lim_{n\to\omega}\tau(x_n)$. Let further $R_\omega= R'\cap R^\omega$, where $R$ is embedded into $R^\omega$ as constant sequences. Recall that 
$R_\omega$, $R\vee R_\omega$ are factors and 
$(R\vee R_\omega)'\cap R^\omega = {\bf C}$. 

Also, if $\teta$ is an automorphism of $R$ then 
there exists a unitary element 
$u_\teta\in R^\omega$ normalizing both 
$R$ and $R\vee R_\omega$ such that $\mathrm{Ad}(u)_{|R}=\teta$. 
The unitary element $u_\teta$ is unique 
modulo perturbation by a unitary from $R_\omega$ 
and $\teta$ is properly outer if and only if 
$\mathrm{Ad}(u)$ is 
properly outer on $R\vee R_\omega$ 
(see Connes~\cite{Connes76} for all this). Thus, if $\sigma$ is a properly outer 
cocycle action of some group $G$ on $R$ then $\{\mathrm{Ad}(u_{\sigma(g)}) 
\mid g \in G\}\subset \mathcal{N}(R\vee R_\omega)$ 
implements a cocycle action $\widetilde{\sigma}$ of $G$ 
on $R\vee R_\omega$ (again, see~\cite{Connes76}).

\begin{theorem}
\ 

\noindent
(i)~$R^\omega$ is amenable both relative to $R_\omega$ and to $R\vee R_\omega$.

\smallskip
\noindent
(ii)~If $N$ is the von Neumann algebra generated by the normalizer 
$\mathcal{N}(R\vee R_\omega)$ of $R\vee R_\omega$ in $R^\omega$, 
then $R\vee R_\omega$ is not co-amenable in $N$. Moreover, if $\sigma$ is a properly outer action of an 
infinite property~(T) group $G$ on $R$ (\emph{e.g.} by Bernoulli shifts) 
then $N_0\deq (R\vee R_\omega)\cup \{u_{\sigma(g)} \mid g\in G\})'' 
\simeq (R\vee R_\omega)\rtimes_{\widetilde{\sigma}}G$ 
is not amenable relative to $R\vee R_\omega$. In fact, $N_0$ has property~(T) relative to $R\vee R_\omega$.
\end{theorem}

\begin{proof} By Remark~\ref{rem_alg} it is sufficient to prove that $R^\omega$ 
is amenable relative to $R_\omega$. We show that in fact 
a much stronger version of condition~(iii) in Proposition~\ref{prop_conditions} holds true, namely: 
$\forall\, \mathcal{U}_0 = \{u_n\}_n \subset \mathcal{U}(R^\omega)$ 
countable, $\exists v\in \mathcal{U}(R^\omega)$ such that 
$f=ve_{R_\omega}v^*$ satisfies $ufu^*=f, \forall u\in \mathcal{U}_0$  
or, equivalently, $v\mathcal{U}_0v^* \subset R_\omega$. 

Let $u_n=(u_n^k)_k$ with $u_n^k\in \mathcal{U}(R), \forall k,n$. 
Let $R=\overline{\cup}_pM_{2^p\times 2^p}({\bf C})$ 
and for each $n$ let $p_n$ be such that $\|E_{M_{2^{p_n}\times 2^{p_n}}(\bf C)}
(u_m^k)- u_m^k\|_2 \leq 1/n, \forall k, m \leq n$. 
Thus, if we 
put $P\deq \Pi_{n\to\omega}M_{2^{p_n}\times 2^{p_n}}({\bf C})$ 
then $\mathcal{U}_0, R \subset P$. In particular, 
$P'\cap R^\omega \subset R_\omega$. For each $n$ let now 
$v_n\subset R$ be a unitary element satisfying 
$v(M_{2^{p_n}\times 2^{p_n}}({\bf C}))v_n^* \subset 
M_{2^{p_n}\times 2^{p_n}}({\bf C})'\cap R$ (this 
is trivially possible, since the latter 
is a type $\mathrm{II}_1$ factor). Thus, if we let $v=(v_n)_n$ 
then $vPv^* \subset P'\cap R^\omega \subset R_\omega$, 
showing that $v\mathcal{U}_0v^* \subset R_\omega$ 
as well.     

The rest of the statement is now trivial, 
because if we denote 
by $\mathcal{G}$ the quotient group $\mathcal{N}(R\vee R_\omega)/ \mathcal{U}(R\vee R_\omega)$ 
then $N$ itself is a cross product of the form 
$(R\vee R_\omega) \rtimes \mathcal{G}$.
By Proposition~\ref{prop_equi} (with $G=\mathcal{G}$ and $H=1$), the co-amenability 
of $R\vee R_\omega$ in $N$ amounts 
to the amenability of the group $\mathcal{G}$. 
But $\mathcal{G}$ contains copies of any countable non-amenable 
group (\emph{e.g.} property~(T) groups). 
\end{proof}

%

%
\end{document}